\documentclass{svproc}

\usepackage{graphicx}
\usepackage{booktabs}
\usepackage{amsmath,amssymb}
\usepackage{multirow}
\usepackage{url}

\usepackage[hidelinks]{hyperref}

\graphicspath{{figures/}}

\begin{document}
\mainmatter

\title{Stochastic Mixed-Integer Optimization of Dynamic Electricity Tariffs with Consumer Protection}
\titlerunning{Consumer-Protected Dynamic Tariff Optimization}

\author{
Ananya Kale\inst{1} \and
Mohit Apte\inst{2} \and
Chhaya Gosavi\inst{1}
}
\authorrunning{Ananya Kale et al.}
\tocauthor{Ananya Kale, Mohit Apte and Chhaya Gosavi}
\institute{
MKSSS's Cummins College of Engineering for Women, India\\
\email{ananya.kale@cumminscollege.in}\
\email{chhaya.gosavi@cumminscollege.in}
\and
University of Chicago, Chicago, IL, USA\\
\email{mohitapte@uchicago.edu}
}

\maketitle

\begin{abstract}
Day-ahead residential tariffs must be posted before demand is observed.
Aggressive high prices can cut peaks in simulation, but they can also raise
customer bills. This paper measures how much peak reduction is lost when a
tariff design problem is required to keep revenue near a flat-tariff baseline
and to limit bill increases. Using half-hourly data from 5{,}566 Low Carbon
London households (167.8 million validated readings), we estimate
quasi-experimental price response against the standard-tariff comparison
group, form bootstrap demand scenarios, and solve stochastic mixed-integer
programs with HiGHS. On all 73 eligible held-out test days, a
segment-protected stochastic tariff reduces simulated peak demand by
2.29\% on average (95\% day-bootstrap CI $[2.11,2.48]$), with revenue change
$-2.26\%$ and mean worst-segment bill increase $2.48\%$. Removing the
segment bill cap raises peak reduction only to $2.37\%$, while the worst
segment bill increase rises to $6.81\%$: in this sample, substantial average
protection costs little peak-shaving performance. The same schedules leave a
household 95th-percentile bill increase of $8.61\%$ (CVaR$_{95}$ $14.76\%$),
so segment-average caps do not bound household tails. We report the full
price-of-protection frontier under consistent household simulation and
compare segment protection with a representative-household (tail-aware)
variant: the latter cuts household p95 from $8.61\%$ to $7.05\%$ while
changing mean peak reduction only from $2.29\%$ to $2.28\%$.
All optimized outcomes are model-based counterfactuals under an
opt-in trial; wholesale costs are unavailable, so we do not optimize profit.
\end{abstract}

\section{Introduction}

A supplier that posts a day-ahead time-of-use (ToU) tariff chooses one of a
small set of price levels for each half-hour before consumption is realized.
The operational hope is that high prices in peak half-hours shift load and
reduce system peaks. The design problem is constrained in practice: revenue
must stay near the level earned under a flat normal tariff, schedules cannot
be arbitrarily complex, and customer bills cannot rise without bound.

These constraints matter. On the Low Carbon London (LCL) trial data used here,
a simple rule that places high prices on forecast peaks reduces simulated peak
demand by about $5.6\%$ relative to a flat normal tariff, but it also raises
mean segment bills by more than $100\%$. A revenue-neutral schedule with a
$+3\%$ segment bill cap achieves a smaller peak reduction---about $2.3\%$---with
bill increases held near $2.5\%$ at the segment mean. The gap between those
two numbers is the quantity this paper studies: the operational price of
consumer protection under demand and response uncertainty.

Two further observations sharpen the claim.
First, removing the segment bill cap changes mean peak reduction only from
$2.29\%$ to $2.37\%$ on the full held-out test period, while the worst-segment
bill increase rises from $2.48\%$ to $6.81\%$. In this dataset, average
consumer protection is cheap in peak-shaving terms.
Second, the same protected schedules leave a household 95th-percentile bill
increase of $8.61\%$ and a household CVaR$_{95}$ of $14.76\%$. Segment-mean
caps therefore do not describe the upper tail of household bill impacts under
heterogeneous load shapes.

The contributions are empirical and operational:
\begin{enumerate}
\item full held-out evaluation of day-ahead tariff MILPs on every eligible LCL
test day (73 days), with paired day-bootstrap intervals;
\item quantification of the mismatch between segment-average bill protection
and household bill tails, including a protection frontier and a
representative-household (tail-aware) protection variant;
\item an empirical price-of-protection frontier under demand and response
uncertainty that links bill caps to peak reduction, revenue, and household
tail risk.
\end{enumerate}
As a robustness check, we also re-estimate response with lead and rebound
terms and embed a state-transition MILP; estimated intertemporal effects are
small, and the protected peak-reduction result is essentially unchanged.
The software pipeline is documented separately and is not claimed as a
research contribution.

\section{Related work}

Empirical work on residential dynamic pricing documents average reductions in
peak-period use under critical-peak and ToU rates, with substantial
heterogeneity across pilots~\cite{faruqui2010,herter2007,allcott2011,jessoe2014}.
Survey and meta-analytic evidence emphasizes that bill impacts and engagement
vary sharply across households even when average peak reductions are
modest~\cite{faruqui2017,joskow2012}.
The Low Carbon London trial provides a large London AMI sample with a
three-level dynamic tariff, day-ahead notification, and a contemporaneous
standard-tariff comparison group~\cite{lcl2014,schofield2014,schofield2015}.
Because LCL recruitment was opt-in, the standard group is not a randomized
control; we treat differentials as quasi-experimental.
Related identification work on nonlinear and average-price perception
\cite{ito2014,wolak2011} motivates caution in extrapolating trial elasticities
outside the observed price ladder.

On the optimization side, day-ahead and ToU design problems are often cast as
mixed-integer programs with peak, cost, or comfort
objectives~\cite{celebi2012,de2013,zugno2013,samadi2010}.
Stochastic and risk-aware formulations appear when demand or prices are
uncertain~\cite{conejo2010,nojavan2021,jamshidi2021}, and Rockafellar--Uryasev
CVaR is a standard linearizable risk measure~\cite{rockafellar2000}.
Retailer and aggregator models commonly protect expected profit or a
representative consumer, rather than the upper tail of household
bills~\cite{hatami2009,garcia2005}.

Fairness and affordability constraints have been studied in retail pricing and
demand-response design~\cite{hupez2021,moret2020,burger2019}.
In practice, however, computational models still often impose aggregate or
segment-average bill caps as a proxy for household protection.
Protection of expected segment bills does not necessarily bound the upper
tail of household-level bill impacts under heterogeneous consumption and
response. We quantify this gap empirically on LCL data and measure the
operational cost of progressively stronger consumer protection under
uncertainty.

\section{Data and empirical models}

\subsection{Data}

We use the partitioned LCL half-hourly release: 168 archive members, 5{,}566
households after ingestion, and 167{,}817{,}021 validated meter rows.
The 2013 dynamic tariff takes values in $\{\mathrm{low},\mathrm{normal},
\mathrm{high}\}$ with prices $3.90$, $11.76$, and $67.20$~p/kWh from the trial
documentation. Weather is Meteostat station data for 2012--2014 with a fixed
filename-to-year mapping validated against file contents.
Households are segmented with $k$-means on pre-treatment 2012 load-shape
features; undersized segments are merged, leaving four segments ($S1$--$S4$)
used in forecasting and optimization.

Chronological splits over 2013 trial days are $60/20/20$: training through
7~August, validation through 19~October, and testing from 20~October through
31~December (73 test days). The 2012 year is used for pre-treatment
segmentation and placebo checks, not for fitting the main response model.

\subsection{Baseline demand and price response}

Segment-level consumption under the normal tariff is forecast with a seasonal
naive benchmark, a regularized linear model, and a histogram gradient-boosting
regressor. Gradient boosting is preferred (validation MAE $0.0099$~kWh per
household). On the test set it attains MAE $0.0140$, RMSE $0.0200$, and WAPE
$6.7\%$.

Let $D_{g,t}$ be the log difference between mean ToU and mean standard-tariff
consumption in segment $g$ at half-hour $t$. We regress $D_{g,t}$ on high/low
indicators interacted with time-of-day blocks and segments, half-hour and month
effects, and weather controls. Standard errors are clustered by date.
Estimation uses training and validation days only. The high-price evening
coefficient is $-0.051$ (SE~$0.0085$); the low-price evening coefficient is
$+0.042$. Directionally consistent estimates appear in a within-ToU predictive
model and in a household two-way fixed-effects check on a stratified sample.
A 2012 placebo that assigns the 2013 tariff calendar to pre-treatment dates is
reported in the repository diagnostics; we do not claim full causal
identification.

\begin{figure}[t]
\centering
\includegraphics[width=0.95\textwidth]{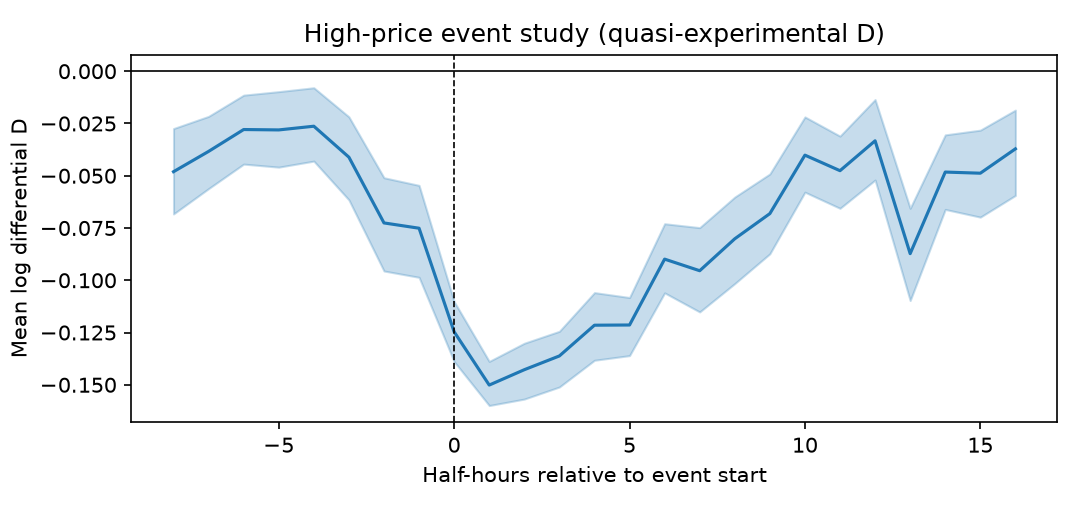}
\caption{Mean log differential $D$ around high-price event starts, with
day-cluster bootstrap bands. Relative half-hour $0$ is the first high-price
interval of an event.}
\label{fig:event}
\end{figure}

For the robustness extension in Section~\ref{sec:robust}, we also estimate lead
and post-high differentials. After removing collinear event indicators,
anticipation and first-period rebound coefficients are small (order
$10^{-3}$). Event-study plots of $D$ around high-price starts are shown in
Fig.~\ref{fig:event}.

\section{Day-ahead tariff optimization}
\label{sec:model}

This section states the stochastic MILP used in the main experiments.
Indices, data, decisions, the objective, and every binding constraint class
are written explicitly so the model can be reconstructed from the text.

\subsection{Indices and data}

For a candidate day let
$t\in\mathcal{T}=\{0,\ldots,47\}$ index half-hours,
$k\in\mathcal{K}=\{\mathrm{low},\mathrm{normal},\mathrm{high}\}$ tariff levels,
$g\in\mathcal{G}$ segments, and
$s\in\mathcal{S}=\{1,\ldots,S\}$ demand scenarios ($S=50$ in the conference
runs).
Prices $p_k$ are in GBP/kWh.
Demand coefficients $d_{g,t,k,s}$ (kWh per household) are precomputed from
baseline forecasts, day-block residual draws, and bootstrap response
multipliers.
Segment weights $w_g$ are ToU household counts.
Define system load, revenue, and segment bills under schedule $x$:
\begin{align}
L_{t,s}(x)
&= \sum_{g\in\mathcal{G}}\sum_{k\in\mathcal{K}} w_g\,d_{g,t,k,s}\,x_{t,k},
\label{eq:load}\\
\mathrm{Rev}_s(x)
&= \sum_{t\in\mathcal{T}}\sum_{g\in\mathcal{G}}\sum_{k\in\mathcal{K}}
   w_g\,d_{g,t,k,s}\,p_k\,x_{t,k},
\label{eq:rev}\\
B_{g,s}(x)
&= \sum_{t\in\mathcal{T}}\sum_{k\in\mathcal{K}}
   d_{g,t,k,s}\,p_k\,x_{t,k}.
\label{eq:bill}
\end{align}
Scenario peak and baseline (all-normal) quantities are
\begin{align}
P_s(x)&=\max_{t\in\mathcal{T}} L_{t,s}(x),\\
P_s^0&=\max_{t\in\mathcal{T}} L_{t,s}(x^{\mathrm{flat}}),\qquad
R^0=\mathbb{E}_s[\mathrm{Rev}_s(x^{\mathrm{flat}})],\qquad
B_g^0=\mathbb{E}_s[B_{g,s}(x^{\mathrm{flat}})],
\end{align}
where $x^{\mathrm{flat}}$ posts normal in every half-hour and
$\mathbb{E}_s$ denotes the equal-weight scenario average.
Let $P^0=\mathbb{E}_s[P_s^0]$.

\subsection{Decision variables}

Binary tariff selection:
\begin{equation}
x_{t,k}\in\{0,1\},\qquad
\sum_{k\in\mathcal{K}} x_{t,k}=1\qquad\forall t\in\mathcal{T}.
\label{eq:onelevel}
\end{equation}
Auxiliary continuous variables:
peak $P_s\ge L_{t,s}(x)$ for all $t,s$;
Rockafellar--Uryasev auxiliaries $\eta\in\mathbb{R}$ and
$u_s\ge 0$ for CVaR of peak;
ramp magnitudes $r_{t,s}\ge |L_{t,s}-L_{t-1,s}|$;
transition linearization $v_{t,k}\ge |x_{t,k}-x_{t-1,k}|$, with
\begin{equation}
T_{\mathrm{tr}}(x)=\tfrac12\sum_{t=1}^{47}\sum_{k\in\mathcal{K}} v_{t,k}.
\label{eq:trans}
\end{equation}
CVaR of peak at level $\alpha=0.9$ is
\begin{equation}
\mathrm{CVaR}_{\alpha}(P)
=\eta+\frac{1}{(1-\alpha)S}\sum_{s\in\mathcal{S}} u_s,\qquad
u_s\ge P_s-\eta.
\label{eq:cvar}
\end{equation}

\subsection{Objective}

The primary model minimizes a normalized combination of expected peak, peak
CVaR, expected ramp, and tariff transitions:
\begin{equation}
\begin{aligned}
\min_{x}\quad
&\frac{\mathbb{E}_s[P_s]}{P^0}
+\lambda_{\mathrm{CVaR}}\frac{\mathrm{CVaR}_{0.9}(P)}{P^0}\\
&\qquad
+\lambda_R\frac{\mathbb{E}_s\bigl[\sum_{t=1}^{47} r_{t,s}\bigr]}{48\,P^0}
+\lambda_T\frac{T_{\mathrm{tr}}(x)}{48},
\end{aligned}
\label{eq:obj}
\end{equation}
with default weights
$(\lambda_{\mathrm{CVaR}},\lambda_R,\lambda_T)=(0.5,0.1,0.05)$.
There is no wholesale procurement cost in the data, so the model does not
optimize profit or gross margin.

\subsection{Consumer-protection and complexity constraints}

Revenue neutrality around flat-normal expected revenue
($\delta=0.03$ by default):
\begin{equation}
(1-\delta)R^0
\le \mathbb{E}_s[\mathrm{Rev}_s(x)]
\le (1+\delta)R^0.
\label{eq:revband}
\end{equation}
Segment expected-bill protection (cap $\gamma$; default $\gamma=0.03$):
\begin{equation}
\mathbb{E}_s[B_{g,s}(x)]
\le (1+\gamma)\,B_g^0
\qquad\forall g\in\mathcal{G}.
\label{eq:billcap}
\end{equation}
Optional scenario-tail bill protection (used in the tail-aware variant), with
auxiliaries $(\eta_g^B,u_{g,s}^B)$ and level $\alpha=0.9$:
\begin{equation}
\mathrm{CVaR}_{\alpha}(B_g)
\le (1+\gamma^{\mathrm{tail}})\,B_g^0
\qquad\forall g\in\mathcal{G}.
\label{eq:billcvar}
\end{equation}
Event-frequency and complexity limits:
\begin{align}
\sum_{t\in\mathcal{T}} x_{t,\mathrm{high}} &\le H_{\max},
&
\sum_{t\in\mathcal{T}} x_{t,\mathrm{low}} &\le L_{\max},
\label{eq:freq}\\
T_{\mathrm{tr}}(x) &\le T_{\max},
\label{eq:maxtrans}\\
\sum_{\tau=t}^{t+C} x_{\tau,\mathrm{high}} &\le C
\qquad\forall t\text{ with }t+C\le 47,
\label{eq:consec}
\end{align}
with defaults $(H_{\max},L_{\max},T_{\max},C)=(12,24,10,6)$.
Minimum run length $L_{\mathrm{run}}=2$ for non-normal levels: if level $k\neq\mathrm{normal}$
turns on at $t\ge 1$, it must remain on for the next $L_{\mathrm{run}}-1$ periods,
\begin{equation}
x_{\tau,k}
\ge x_{t,k}-x_{t-1,k}
\qquad
\forall\,k\neq\mathrm{normal},\;
\tau=t+1,\ldots,t+L_{\mathrm{run}}-1.
\label{eq:minrun}
\end{equation}

\subsection{Intertemporal state-transition extension}
\label{sec:itmath}

When load depends on the previous and current tariff levels, introduce
transition binaries $z_{t,k^-,k}\in\{0,1\}$ with
\begin{align}
\sum_{k^-,k} z_{t,k^-,k}&=1,
&
x_{t,k}&=\sum_{k^-} z_{t,k^-,k},
\label{eq:zlink}\\
\sum_{k^-} z_{t,k^-,k}
&=\sum_{k^+} z_{t+1,k,k^+},
&
\sum_{k} z_{0,\mathrm{normal},k}&=1.
\label{eq:zchain}
\end{align}
Precomputed coefficients $d_{g,t,k^-,k,s}$ absorb contemporaneous multipliers and
clipped lead/rebound adjustments. System load and revenue are then linear in
$z$; bill and revenue constraints continue to use the same-period coefficients
in~\eqref{eq:bill}--\eqref{eq:revband} as a conservative approximation.
The resulting MILP remains solvable with open-source HiGHS (mean solve time
about $19$~s versus about $7$~s for the same-period model).

\subsection{Tail-aware protection variant}

Segment constraints protect averages. To probe household tails inside the
optimizer, we select $A=8$ representative ToU households (highest evening load,
balanced across segments), append them as zero-weight archetypes so that they
enter bill constraints but not system peak or revenue, and impose expected-bill
caps of $\gamma=0.03$ on segments together with a looser archetype cap
$\gamma^{\mathrm{arch}}=0.10$ and scenario CVaR bill caps at
$\gamma^{\mathrm{arch}}$ on archetype indices only.
Enforcing the same $+3\%$ mean cap on peaky archetypes collapses simulated peak
reduction toward zero; that extreme case underscores the cost of literal
household-style protection.
Peak, revenue, and segment KPIs are scored on the original segment data after
optimization so comparisons remain commensurate.

\section{Experimental design}

Every chronological test day with complete panel, weather, and lag features is
eligible; all 73 test days qualify. Policies scored on the same days and
scenarios include: flat normal; historical trial tariff; rule-based
high-on-peak / low-on-trough; deterministic and stochastic peak minimization
with segment protection; stochastic optimization without bill caps; robust
(worst-case) peak minimization; and the tail-aware variant above.
Uncertainty scenarios combine validation residual blocks with bootstrap draws
of response coefficients. Inference uses the day as the sampling unit: paired
differences versus flat normal are summarized by means and day-bootstrap 95\%
intervals. Optimized schedules that were never observed in the trial are
labeled as simulated throughout.

The intellectual progression in the results is
\emph{unprotected pricing $\rightarrow$ segment protection $\rightarrow$
household-tail evidence (and tail-aware protection)}.
Separately, we report deterministic versus stochastic versus robust models as
an uncertainty experiment, and the intertemporal formulation as a behavioral
robustness check.

The protection frontier varies the segment bill cap over
$\{\infty,5\%,3\%,2\%,1\%,0.5\%\}$ on all eligible test days. Household bill
tails at each frontier point are simulated with a fixed household sample and
day-averaged household percentage changes, so percentiles are comparable
across caps. Policy stability is measured by redrawing response coefficients
and reoptimizing on five representative days (25 draws each).

\section{Results}

\subsection{Unprotected pricing versus segment protection}

\begin{figure}[t]
\centering
\includegraphics[width=0.95\textwidth]{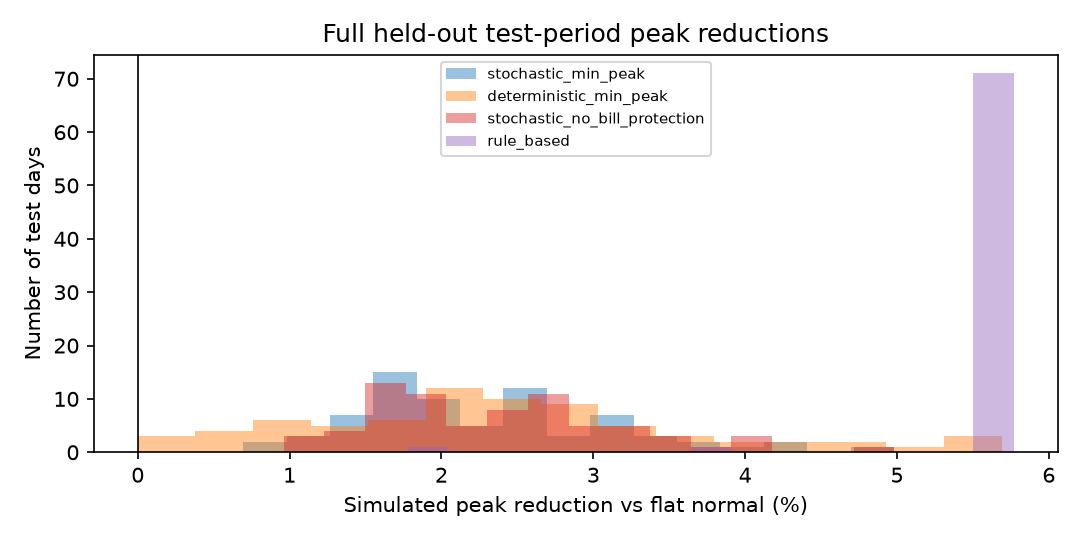}
\caption{Distribution across 73 test days of simulated peak reduction versus
flat normal for selected policies.}
\label{fig:peakdist}
\end{figure}

Table~\ref{tab:policy} reports mean simulated peak reduction, revenue change,
and worst-segment bill increase over 73 test days.
Figure~\ref{fig:peakdist} shows the day-level distribution of peak reductions
for selected policies.

\begin{table}[t]
\caption{Simulated out-of-sample policy comparison on 73 held-out test days.
Peak and revenue changes are relative to flat normal. Intervals are
day-bootstrap 95\% CIs for mean peak reduction. Outcomes for optimized and
rule-based schedules are model-based.}
\label{tab:policy}
\centering
\scriptsize
\begin{tabular}{@{}l r r r r@{}}
\toprule
Policy & Peak red.\ (\%) & 95\% CI & Rev.\ $\Delta$ (\%) & Max seg.\ bill (\%) \\
\midrule
Flat normal & 0.00 & --- & 0.00 & 0.00 \\
Historical dToU & 0.26 & $[-0.22,0.76]$ & 20.16 & 24.73 \\
Rule-based & 5.61 & $[5.50,5.69]$ & 138.19 & 159.52 \\
Deterministic min-peak & 2.39 & $[2.12,2.69]$ & $-2.06$ & 2.09 \\
Stochastic, segment cap & 2.29 & $[2.11,2.48]$ & $-2.26$ & 2.48 \\
Stochastic, no bill cap & 2.37 & $[2.20,2.54]$ & 1.20 & 6.81 \\
Tail-aware (arch.\ $+10\%$) & 2.28 & $[2.09,2.47]$ & $-2.56$ & 2.44 \\
Robust worst-case & 1.37 & $[1.12,1.64]$ & $-2.08$ & 2.11 \\
Intertemporal (robustness) & 2.55 & $[2.38,2.73]$ & $-2.25$ & 2.50 \\
\bottomrule
\end{tabular}
\end{table}

\begin{figure}[t]
\centering
\includegraphics[width=0.95\textwidth]{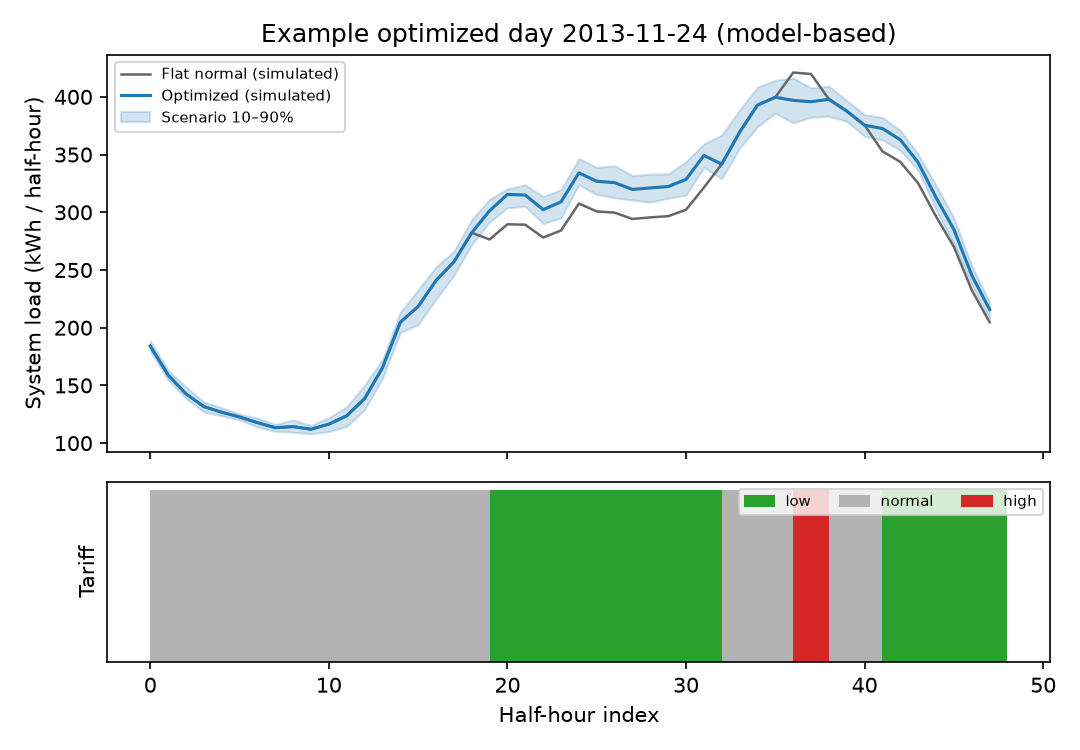}
\caption{Example held-out day: simulated system load under flat normal and
under the segment-protected stochastic schedule, with scenario band and
posted tariff levels.}
\label{fig:example}
\end{figure}

The headline comparison is between the segment-protected stochastic policy and
the same model without a bill cap. Removing protection raises mean peak
reduction only from $2.29\%$ to $2.37\%$---a gain of $0.08$ percentage
points---while the worst-segment bill increase rises from $2.48\%$ to
$6.81\%$. In this sample, substantial average consumer protection can be
obtained for little sacrifice in peak-shaving performance.
Rule-based high-on-peak pricing delivers a larger simulated peak cut
($5.61\%$) but is not a credible consumer-facing schedule: revenue and segment
bills inflate by more than $100\%$. An example optimized day is shown in
Fig.~\ref{fig:example}.

\subsection{Segment protection versus household bill tails}

Segment constraints bind the optimizer, but they do not describe household
tails. Under the $+3\%$ segment-protected stochastic schedules, simulated
household bill changes (1{,}112 ToU households; day-averaged per household)
have overall p90 $5.25\%$, p95 $8.61\%$, p99 $19.43\%$, and CVaR$_{95}$
$14.76\%$ (Table~\ref{tab:hh}).
About $14.7\%$ of households exceed $+3\%$, $10.6\%$ exceed $+5\%$, and
$3.5\%$ exceed $+10\%$.
Segment $S2$ is the worst group: p95 $12.86\%$ and CVaR$_{95}$ $21.20\%$, with
$30.8\%$ of households above $+3\%$.

\begin{table}[t]
\caption{Household bill-risk under the segment-protected stochastic policy
(day-averaged bill change vs flat normal; $n=1{,}112$ households).}
\label{tab:hh}
\centering
\scriptsize
\begin{tabular}{@{}l r r r r r r r@{}}
\toprule
Group & p90 & p95 & p99 & CVaR$_{95}$ & $>3\%$ & $>5\%$ & $>10\%$ \\
\midrule
Overall & 5.25 & 8.61 & 19.43 & 14.76 & 0.147 & 0.106 & 0.035 \\
$S1$ & 2.53 & 6.71 & 13.17 & 12.53 & 0.085 & 0.071 & 0.023 \\
$S2$ & 8.80 & 12.86 & 24.21 & 21.20 & 0.308 & 0.219 & 0.086 \\
$S3$ & 3.52 & 5.94 & 14.53 & 10.97 & 0.115 & 0.066 & 0.025 \\
$S4$ & 2.52 & 5.87 & 9.41 & 8.91 & 0.084 & 0.061 & 0.009 \\
\bottomrule
\end{tabular}
\end{table}

\begin{figure}[t]
\centering
\includegraphics[width=0.85\textwidth]{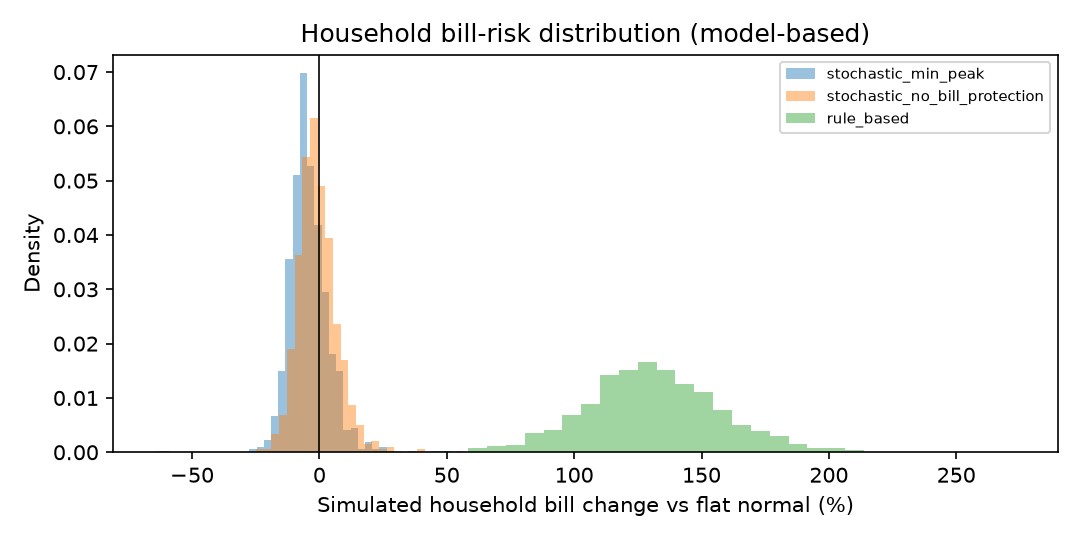}
\caption{Distribution of simulated household bill changes versus flat normal
under selected policies. Segment-protected stochastic schedules still leave a
right tail of bill increases.}
\label{fig:hhbills}
\end{figure}

The gap between a $3\%$ segment cap and an $8.61\%$ household p95 is the
empirical content of the paper's title: average protection is not tail
protection.

\subsection{Price of protection}

\begin{figure}[t]
\centering
\includegraphics[width=0.85\textwidth]{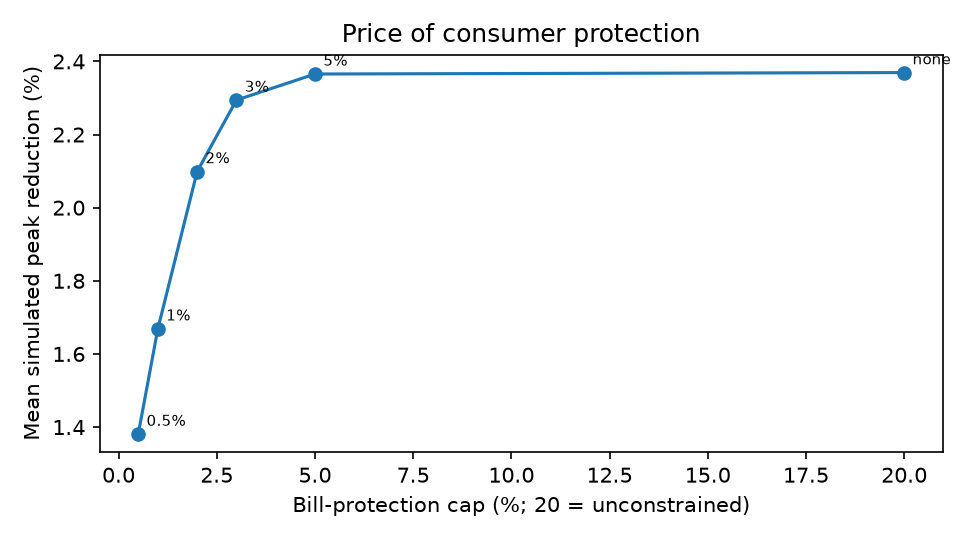}
\caption{Mean simulated peak reduction as the segment bill-protection cap
tightens (``20'' on the horizontal axis marks the unconstrained case).}
\label{fig:protection}
\end{figure}

Figure~\ref{fig:protection} and Table~\ref{tab:frontier} summarize the
protection frontier under consistent household simulation on all 73 test days.
Tightening the segment bill cap from unconstrained to $0.5\%$ reduces mean
simulated peak reduction from $2.37\%$ to $1.38\%$.
Household percentiles move with the segment cap once measurement is held fixed:
household p95 falls from $12.22\%$ (unconstrained) to $8.61\%$ under a $3\%$
segment cap and to $6.45\%$ under a $0.5\%$ cap.
Earlier stride-subset runs that mixed full-sample and pooled day-level
percentiles produced non-monotonic household tails; those artifacts disappear
under the uniform day-averaged household design used here.
The remaining gap is still large: a $3\%$ segment-mean cap coexists with an
$8.61\%$ household p95 and a $14.76\%$ household CVaR$_{95}$.

Adding representative-household constraints (eight high-evening-load archetypes
with $+10\%$ expected-bill and CVaR caps, on top of the $+3\%$ segment mean
cap) changes mean peak reduction only from $2.29\%$ to $2.28\%$, while
household p95 falls from $8.61\%$ to $7.05\%$ and CVaR$_{95}$ from $14.76\%$
to $13.61\%$ (Table~\ref{tab:hhcompare}).
The share of households above $+5\%$ falls from $10.6\%$ to $8.4\%$.
In this sample, moving from segment-average protection toward
household-tail-aware protection is cheap in peak-shaving terms.
Enforcing the same $+3\%$ mean cap on those archetypes instead of $+10\%$
collapses peak reduction toward zero in pilot solves: literal household-style
caps at the segment rate are much more expensive.

\begin{table}[t]
\caption{Segment protection versus tail-aware (archetype) protection on 73
test days. Household tails are day-averaged simulated bill changes
($n=1{,}112$).}
\label{tab:hhcompare}
\centering
\scriptsize
\begin{tabular}{@{}l r r r r r@{}}
\toprule
Policy & Peak red.\ (\%) & HH p95 (\%) & HH CVaR$_{95}$ (\%) & HH $>5\%$ & HH $>10\%$ \\
\midrule
No bill cap & 2.37 & 12.22 & 18.58 & 0.190 & 0.082 \\
Segment cap $+3\%$ & 2.29 & 8.61 & 14.76 & 0.106 & 0.035 \\
Tail-aware (arch.\ $+10\%$) & 2.28 & 7.05 & 13.61 & 0.084 & 0.030 \\
\bottomrule
\end{tabular}
\end{table}

\begin{table}[t]
\caption{Protection frontier on all 73 eligible held-out test days. Household
percentiles are day-averaged simulated bill increases under the optimized
schedules (fixed sample of 1{,}112 households; consistent across caps).}
\label{tab:frontier}
\centering
\scriptsize
\begin{tabular}{@{}l r r r r r@{}}
\toprule
Bill cap & Peak red.\ (\%) & Rev.\ $\Delta$ (\%) & HH p95 (\%) & HH CVaR$_{95}$ (\%) & HH $>10\%$ \\
\midrule
None & 2.37 & 1.20 & 12.22 & 18.58 & 0.082 \\
$5\%$ & 2.37 & $-1.03$ & 9.87 & 16.04 & 0.048 \\
$3\%$ & 2.29 & $-2.26$ & 8.61 & 14.76 & 0.035 \\
$2\%$ & 2.10 & $-2.63$ & 7.84 & 13.97 & 0.032 \\
$1\%$ & 1.67 & $-2.56$ & 6.66 & 12.22 & 0.025 \\
$0.5\%$ & 1.38 & $-2.50$ & 6.45 & 11.79 & 0.024 \\
\bottomrule
\end{tabular}
\end{table}

\begin{figure}[t]
\centering
\includegraphics[width=0.85\textwidth]{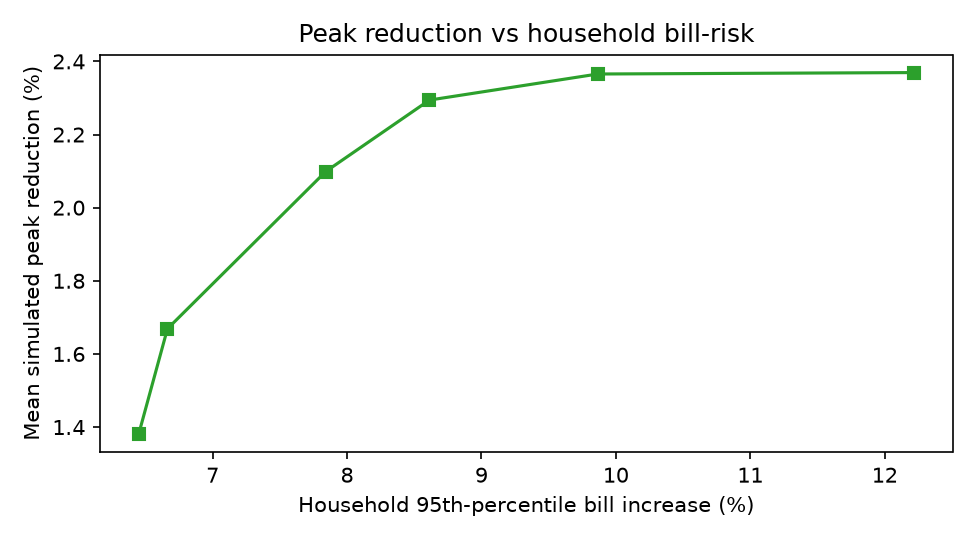}
\caption{Peak reduction versus household 95th-percentile bill increase along
the protection frontier.}
\label{fig:peakhh}
\end{figure}

\subsection{Uncertainty and behavioral robustness}
\label{sec:robust}

\begin{figure}[t]
\centering
\includegraphics[width=0.95\textwidth]{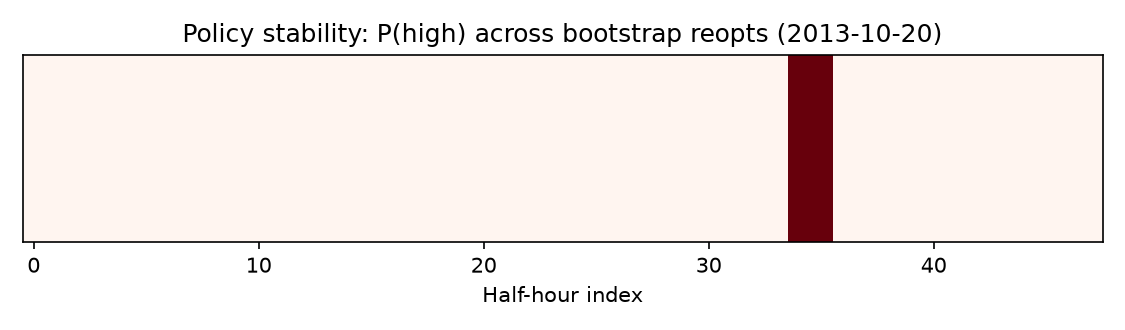}
\caption{Share of bootstrap reoptimizations that select a high price in each
half-hour on a representative test day.}
\label{fig:stability}
\end{figure}

Deterministic optimization achieves a slightly higher mean peak reduction
($2.39\%$) than the stochastic model ($2.29\%$). The robust worst-case model
is more conservative ($1.37\%$). Bootstrap reoptimization agreement across
tariff half-hours averages $0.82$ on the sampled days
(Fig.~\ref{fig:stability}).

With the intertemporal state-transition model of
Section~\ref{sec:itmath}, mean peak reduction is $2.55\%$
(CI $[2.38,2.73]$), slightly above the same-period stochastic policy.
Evaluating those schedules under same-period demand yields $2.25\%$, so the
schedule change matters, but estimated lead and first-period rebound
coefficients remain small (Table~\ref{tab:it}).
The useful reading is robustness: the main protection result survives a richer
behavioral specification, and a same-period model is already a reasonable
approximation for this dataset.

\begin{table}[t]
\caption{Selected intertemporal difference-regression coefficients
(date-clustered SEs). Same-period evening high/low references are shown for
comparison.}
\label{tab:it}
\centering
\scriptsize
\begin{tabular}{@{}l r r@{}}
\toprule
Term & Estimate & SE \\
\midrule
\texttt{high\_f1} (lead 1) & $-0.0021$ & $0.0033$ \\
\texttt{high\_f2} (lead 2) & $-0.0410$ & $0.0065$ \\
\texttt{post\_high\_1} & $-0.0010$ & $0.0062$ \\
\texttt{post\_high\_2} & $-0.0097$ & $0.0075$ \\
\texttt{post\_high\_4} & $+0.0118$ & $0.0075$ \\
\texttt{high\_x\_evening} (IT) & $-0.0175$ & $0.0118$ \\
\texttt{high\_x\_evening} (same-period) & $-0.0511$ & $0.0085$ \\
\bottomrule
\end{tabular}
\end{table}

\section{Limitations}

LCL ToU enrollment was opt-in; the standard-tariff group is a comparison group,
not a randomized control. Response is identified only at three observed price
levels. Optimized outcomes are simulated counterfactuals, not field experiments.
Optimization uses realized weather and is therefore an ex-post information
benchmark. No wholesale prices are available, so we report tariff revenue and
load metrics rather than profit. External validity is limited to this London
sample and period. Socioeconomic ACORN labels are absent from the processed
household table used here. Representative-household constraints tighten
protection for selected archetypes; they are not a full household-level
robust optimization over all meters.

\section{Conclusion}

On the full LCL held-out test period, segment-protected stochastic tariffs
produce simulated peak reductions near $2.3\%$, with revenue held within a few
percent of the flat baseline.
Removing segment bill protection buys almost no additional peak reduction while
materially raising segment bills: average protection is inexpensive in
peak-shaving terms here.
The same schedules leave household bill tails several times larger than the
segment cap. A representative-household (tail-aware) extension narrows those
tails at almost no additional peak cost in this sample. The useful planning
object is therefore the feasible tradeoff among peak, revenue, segment bills,
and household bill risk---not unconstrained peak minimization, and not
segment averages alone.



\begin{thebibliography}{99}

\bibitem{allcott2011}
Allcott, H.:
Rethinking real-time electricity pricing.
Resource and Energy Economics 33(4), 820--842 (2011)

\bibitem{burger2019}
Burger, S., Knittel, C., P{\'e}rez-Arriaga, I., Schneider, I., vom~Scheidt, F.:
The efficiency and distributional effects of alternative residential
electricity rate designs.
NBER Working Paper 25570 (2019)

\bibitem{celebi2012}
Celebi, E., Fuller, J.D.:
Time-of-use pricing in electricity markets under different market structures.
IEEE Transactions on Power Systems 27(3), 1170--1181 (2012)

\bibitem{conejo2010}
Conejo, A.J., Carri{\'o}n, M., Morales, J.M.:
Decision Making Under Uncertainty in Electricity Markets.
Springer, New York (2010)

\bibitem{de2013}
de~S{\'a} Ferreira, R., Barroso, L.A., Lino, P.R., Carvalho, M.M.,
Valenzuela, P.:
Time-of-use tariff design under uncertainty in price-elasticities of
electricity demand: a stochastic optimization approach.
IEEE Transactions on Smart Grid 4(4), 2285--2295 (2013)

\bibitem{faruqui2010}
Faruqui, A., Sergici, S.:
Household response to dynamic pricing of electricity: a survey of 15
experiments.
Journal of Regulatory Economics 38(2), 193--225 (2010)

\bibitem{faruqui2017}
Faruqui, A., Sergici, S., Warner, C.:
Arcturus 2.0: A meta-analysis of time-varying rates for electricity.
The Electricity Journal 30(10), 64--72 (2017)

\bibitem{garcia2005}
Garc{\'i}a-Bertrand, R., Conejo, A.J., Gabriel, S.:
Electricity cost minimization under uncertainty.
IEEE Transactions on Power Systems 20(2), 743--752 (2005)

\bibitem{hatami2009}
Hatami, A.R., Seifi, H., Sheikh-El-Eslami, M.K.:
Hedging risks with interruptible load programs for a load serving entity.
Decision Support Systems 48(1), 150--157 (2009)

\bibitem{herter2007}
Herter, K.:
Residential implementation of critical-peak pricing of electricity.
Energy Policy 35(4), 2121--2130 (2007)

\bibitem{hupez2021}
Hupez, M., Toubeau, J.-F., Atzeni, I., De~Gr{\`e}ve, Z., Vall{\'e}e, F.:
Pricing electricity consumption in a fair manner: a bilevel approach.
IEEE Transactions on Power Systems 36(3), 2202--2212 (2021)

\bibitem{ito2014}
Ito, K.:
Do consumers respond to marginal or average price? Evidence from nonlinear
electricity pricing.
American Economic Review 104(2), 537--563 (2014)

\bibitem{jamshidi2021}
Jamshidi~Monfared, H., Ghasemi, A., Loni, A., Marzband, M.:
Optimal stochastic conditional value at risk-based management of a demand
response aggregator considering load uncertainty.
In: 2021 IEEE International Conference on Environment and Electrical
Engineering. IEEE (2021)

\bibitem{jessoe2014}
Jessoe, K., Rapson, D.:
Knowledge is (less) power: experimental evidence from residential energy use.
American Economic Review 104(4), 1417--1438 (2014)

\bibitem{joskow2012}
Joskow, P.L., Wolfram, C.D.:
Dynamic pricing of electricity.
American Economic Review 102(3), 381--385 (2012)

\bibitem{lcl2014}
UK Power Networks:
Low Carbon London Project: final / closedown materials
(including data user guide).
UK Power Networks, London (2014--2015)

\bibitem{moret2020}
Moret, F., Pinson, P.:
Energy collectives: a community and fairness based approach to future
electricity markets.
IEEE Transactions on Power Systems 34(5), 3994--4004 (2019)

\bibitem{nojavan2021}
Nojavan, S., Zare, K., Mohammadi-Ivatloo, B.:
CVaR-based retail electricity pricing in day-ahead scheduling of microgrids
under uncertainty.
Energy (2021)

\bibitem{rockafellar2000}
Rockafellar, R.T., Uryasev, S.:
Optimization of conditional value-at-risk.
Journal of Risk 2(3), 21--41 (2000)

\bibitem{samadi2010}
Samadi, P., Mohsenian-Rad, A.H., Schober, R., Wong, V.W.S., Jatskevich, J.:
Optimal real-time pricing algorithm based on utility maximization for smart
grid.
In: First IEEE International Conference on Smart Grid Communications,
pp.\ 415--420. IEEE (2010)

\bibitem{schofield2014}
Schofield, J.R., Carmichael, R., Tindemans, S., Woolf, M., Bilton, M.,
Strbac, G.:
Residential consumer responsiveness to time-varying pricing.
Report A7 for the ``Low Carbon London'' LCNF project.
Imperial College London (2014)

\bibitem{schofield2015}
Schofield, J., Carmichael, R., Tindemans, S., Woolf, M., Bilton, M.,
Strbac, G.:
Experimental validation of residential consumer responsiveness to dynamic
time-of-use pricing.
In: 23rd International Conference on Electricity Distribution (CIRED) (2015)

\bibitem{wolak2011}
Wolak, F.A.:
Do residential customers respond to hourly prices? Evidence from a dynamic
pricing experiment.
American Economic Review 101(3), 83--87 (2011)

\bibitem{zugno2013}
Zugno, M., Morales, J.M., Pinson, P., Madsen, H.:
A bilevel model for electricity retailers' participation in a demand response
market environment.
Energy Economics 36, 182--197 (2013)

\end{thebibliography}
\end{document}